\documentclass[10pt,english,oneside,twocolumn,a4paper]{article}
\usepackage{graphicx}
\usepackage[utf8]{inputenc}
\usepackage{mathptmx}
\usepackage{ragged2e}
\usepackage[singlelinecheck=false]{caption}
\usepackage[T1]{fontenc}
\usepackage[numbers]{natbib}
\usepackage{amsmath}
\RequirePackage[blocks]{authblk}
\usepackage[english]{babel}
\usepackage{enumitem}
\usepackage{eurosym}
\usepackage{microtype}
\usepackage{natbib}
\usepackage{amsmath,mathtools}
\usepackage{amsfonts}
\usepackage{amssymb}
\usepackage{textcomp}
\usepackage{subfig}
\usepackage{fixltx2e}
\usepackage{xcolor}
\usepackage{url}
\usepackage{xspace}
\usepackage{accents}
\usepackage{ifthen}
\usepackage{bm}
\DeclareMathAlphabet{\mathcal}{OMS}{cmsy}{m}{n}

\makeatletter
\let\ps@plain\ps@empty
\def\@xivpt{14bp}

\def\@sect#1#2#3#4#5#6[#7]#8{%
  \ifnum #2>\c@secnumdepth
    \let\@svsec\@empty
  \else
    \refstepcounter{#1}%
    \protected@edef\@svsec{%
      \ifnum #2<4
        \hb@xt@10mm{\csname the#1\endcsname}\relax
      \else
        \hb@xt@12mm{\csname the#1\endcsname}\relax
      \fi}%
  \fi
  \@tempskipa #5\relax
  \ifdim \@tempskipa>\z@
    \begingroup
      #6{%
        \@hangfrom{\hskip #3\relax\@svsec}%
          \interlinepenalty \@M #8\@@par}%
    \endgroup
    \csname #1mark\endcsname{#7}%
    \addcontentsline{toc}{#1}{%
      \ifnum #2>\c@secnumdepth \else
        \protect\numberline{\csname the#1\endcsname}%
      \fi
      #7}%
  \else
    \def\@svsechd{%
      #6{\hskip #3\relax
      \@svsec #8}%
      \csname #1mark\endcsname{#7}%
      \addcontentsline{toc}{#1}{%
        \ifnum #2>\c@secnumdepth \else
          \protect\numberline{\csname the#1\endcsname}%
        \fi
        #7}}%
  \fi
  \@xsect{#5}}
\renewcommand\LARGE{\@setfontsize\LARGE{16}{20}}
\def\abstract#1{\def\@abstract{#1}}
\def\abstractEn#1{\def\@abstractEn{#1}}
\def\titleEn#1{\def\@titleEn{#1}}
\def\@maketitle{%
  \newpage
  \null
  \let \footnote \thanks
    {\LARGE\bfseries\RaggedRight \@titleEn \par}%
    \vskip 1\baselineskip%
    {\normalsize
      \@author\par}%
    \vskip 1.5\baselineskip%
    {\section*{Abstract}
      \@abstractEn}%
  \par
  \vskip 2.5\baselineskip}

\renewcommand\section{\@startsection {section}{1}{\z@}%
                                   {-3.5ex \@plus -1ex \@minus -.2ex}%
                                   {\baselineskip}%
                                   {\normalfont\Large\bfseries\RaggedRight}}
\renewcommand\subsection{\@startsection{subsection}{2}{\z@}%
                                     {\baselineskip}%
                                     {1ex}%
                                     {\normalfont\large\bfseries\RaggedRight}}
\renewcommand\subsubsection{\@startsection{subsubsection}{3}{\z@}%
                                     {1\baselineskip}%
                                     {3bp}%
                                     {\normalfont\normalsize\bfseries\RaggedRight}}
\renewcommand\paragraph{\@startsection{paragraph}{4}{\z@}%
                                    {1\baselineskip\@plus1ex \@minus.2ex}%
                                    {3bp}%
                                    {\normalfont\normalsize\RaggedRight}}
\renewcommand\subparagraph{\@startsection{subparagraph}{5}{\parindent}%
                                       {3.25ex \@plus1ex \@minus .2ex}%
                                       {-1em}%
                                      {\normalfont\normalsize\bfseries\RaggedRight}}
\parindent\p@
\makeatother
\DeclareCaptionLabelSeparator{enskip}{\enskip}
\titleEn{Reducing the Need for New Lines in Germany's Energy Transition:\\ The Hybrid Transmission Grid Architecture}
\author[1]{Matthias Hotz}
\author[2]{Irina Boiarchuk}
\author[2]{Dominic Hewes}
\author[2]{Rolf Witzmann}
\author[1]{Wolfgang Utschick}
\affil[1]{\,Professur f\"ur Methoden der Signalverarbeitung, Technische Universit\"at M\"unchen, Germany (matthias.hotz@tum.de)}
\affil[2]{\,Professur f\"ur Elektrische Energieversorgungsnetze, Technische Universit\"at M\"unchen, Germany}

\abstractEn{The energy transition will lead to an imbalance in electric power generation and demand along the north-south axis of Germany. To manage this imbalance, the transmission system operators proposed a network development plan that requires several thousand kilometers of new lines, which received extensive opposition. In our recent work, we proposed a landscape-preserving network development approach, i.e., the \emph{hybrid architecture}, which relies on the conversion of some \emph{existing} AC lines and transformers to HVDC operation. In this work, we show that this approach can meet the projected capacity and performance requirements with substantially fewer new lines. Due to the high cost of HVDC technology, the investment volume exceeds the current network development plan, but the preservation of landscape and support of public acceptance may justify that cost premium.}

\newcommand{\fBranchRating}{R_k}
\newcommand{\fBranchUtilization}{U_k}
\newcommand{\fBranchKktMultiplier}{\mu_k}
\newcommand{\fTargetRating}{\tilde{R}_k}

\begin{document}

\maketitle

\section{Introduction}
	\label{sec:introduction}

\begin{figure*}[!t]%
\centering%
\hfill%
\subfloat[German transmission grid with the NEP of TSOs/BNetzA.]{%
\includegraphics[width=0.89\columnwidth]{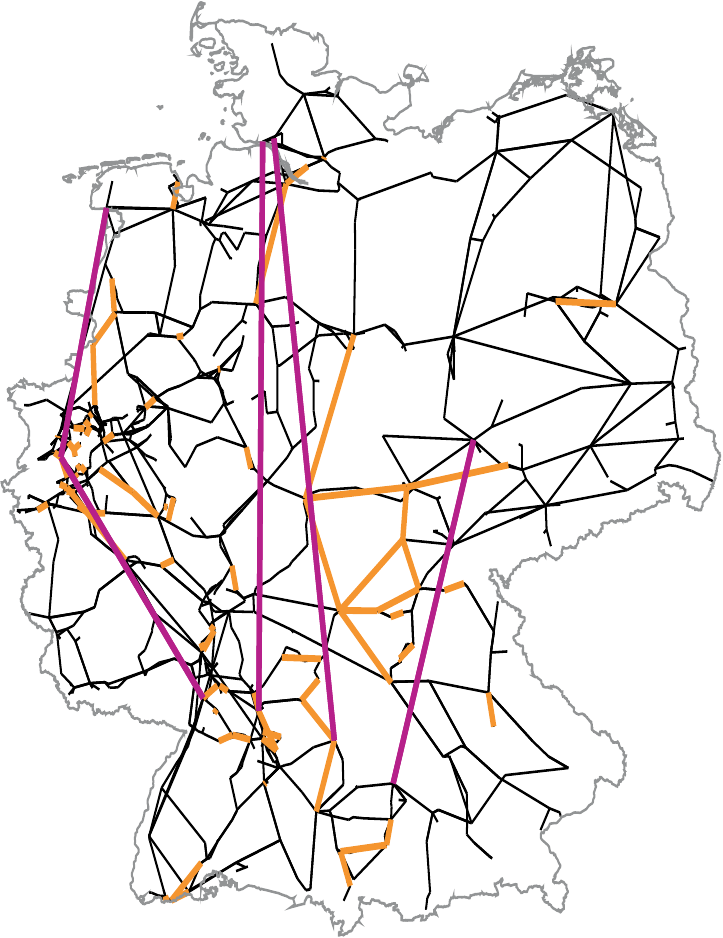}%
\label{fig:grid:tso}}%
\hfill\hfill\hfill%
\subfloat[German transmission grid with the hybrid architecture.]{%
\includegraphics[width=0.89\columnwidth]{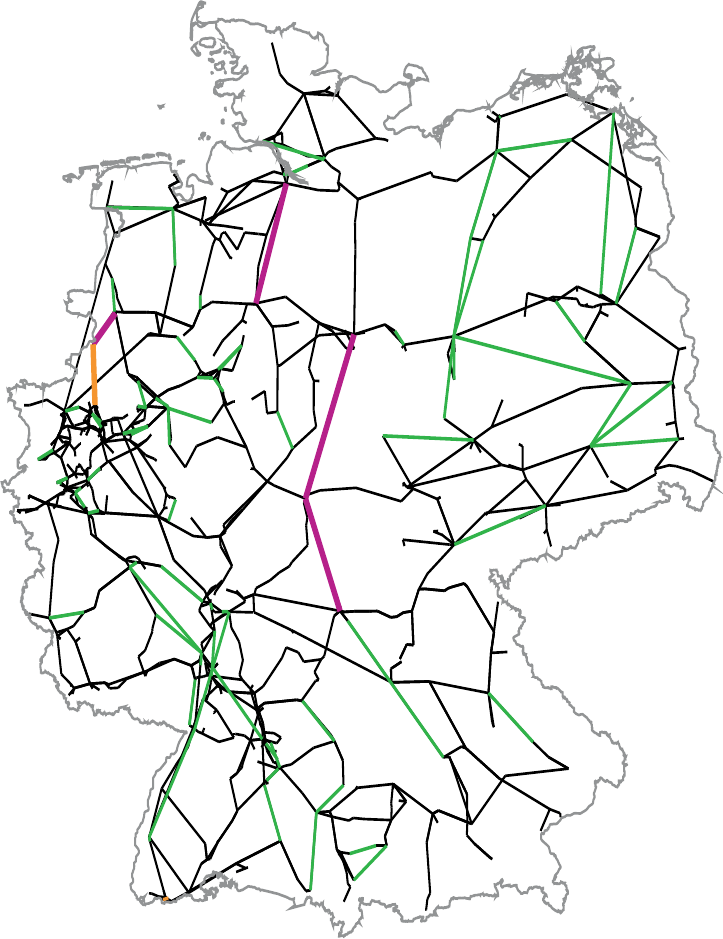}%
\label{fig:grid:htg}}%
\hfill\mbox{}%
\caption{Illustration of the network development strategy of (a) the TSOs/BNetzA as well as (b) the hybrid architecture. Black lines illustrate AC lines, while highlighted lines show AC lines in 2015 that were converted to DC operation in 2030 (green), new DC lines in 2030 (magenta), and new AC lines in 2030 (orange) as identified by step~2) in Fig.~\ref{fig:basegridproc}. In~(b), only $6$ of the $93$ new lines in (a) are retained, reducing the total length of new lines from $5492$\,km to $505$\,km.}%
\label{fig:grid}%
\end{figure*}

The energy transition in Germany is a reorientation of the energy policy towards renewable energy sources (RES) and the reduction of greenhouse gas emissions. In the electricity sector, this reorientation introduces an increasingly distributed and fluctuating energy production due to the growing use of wind, sun, and biomass as energy sources. These changes in the energy mix also entail a geographic shift of power generation, inducing a growing energy surplus in the north and an energy demand in the south of the country~\cite{Netzentwicklungsplan2017a}. This transformation of generation structure renders the expansion of the German transmission grid a key issue of the energy transition~\cite{Netzentwicklungsplan2017a,Bundesnetzagentur2015a}.
Adequate expansion measures are determined via a multi-stage process, which is repeated iteratively since the year 2012~\cite{Bundesnetzagentur2015a}. In this process, the transmission system operators (TSOs) propose projected future scenarios and corresponding network development plans, which are screened, verified, and confirmed by the Federal Network Agency (BNetzA) before specific measures are planned in subsequent stages. Throughout, the proposed measures are repeatedly subject to consultation and public debate to ensure validity and acceptance.

The current network development plan (NEP) is based on projections for the year 2030~\cite{Netzentwicklungsplan2017a}. As in the preceding iterations, it identifies five high-voltage direct current (HVDC) transmission lines as a necessary countermeasure for the north-south generation imbalance. While these HVDC lines are virtually the transmission backbone of the energy transition, they also constitute the primary subject of objection. During consultation of the last confirmed network development plan, more than $34\,000$ statements were received from ministries, agencies, municipalities, associations, organizations, and citizens~\cite{Bundesnetzagentur2015a}. Their statistical evaluation shows that the major concerns comprise the impact on humans, the capital loss of real estate, and the impact on landscape~\cite{Bundesnetzagentur2015a}. These objections may be attributed primarily to the implementation of new transmission corridors, where the majority is due to four of the north-south HVDC lines. Even though the expansion planning follows the so-called NOVA principle~\cite{Netzentwicklungsplan2017a}, i.e., grid optimization and reinforcement is preferred to additional transmission lines, these HVDC lines are identified as necessary measures.

In our previous work, we explored alternative approaches to capacity expansion that rely on \emph{existing} corridors, which resulted in the concept of the \emph{hybrid architecture} introduced in~\cite{Hotz2016a}. The hybrid architecture also takes advantage of HVDC technology, but adopts a different perspective in its utilization. While the current NEP uses \emph{additional} HVDC lines as a \emph{focused} measure against the north-south imbalance, the hybrid architecture comprises a systematic and \emph{system-wide} conversion of certain \emph{existing} AC lines to HVDC. By conversion to HVDC, the transmission capacity of a corridor can be increased by a factor of two or more~\cite{Cigre2014a}. Additionally, voltage source converter (VSC) based HVDC systems offer a rapidly controllable power flow and reactive power capacity at the terminals. Due to this, the system-wide conversion of certain lines does not only selectively increase line capacity, but also introduces substantial \emph{flexibility} to the grid. Furthermore, we have shown that the hybrid architecture induces a transition in \emph{system structure} that supports an efficient utilization~\cite{Hotz2016a,Hotz2017a}. Its structural properties enable a shift of the nonconvex AC optimal power flow (OPF) problem into the \emph{convex} domain, rendering it amenable to efficient solution algorithms and a powerful framework of mathematical theory. Therewith, the hybrid architecture improves decision making, which further supports the efficient utilization of the grid.

In this work, we show that the hybrid architecture is a potential alternative to the NEP. It induces a grid expansion whose performance is on par with the NEP, whilst substantially reducing the need for new lines. Hereafter, Section~\ref{sec:nep} discusses the NEP. Section~\ref{sec:transition} reviews the hybrid architecture and, subsequently, utilizes this concept to develop an alternative network development strategy, whose operational performance is then compared to the NEP. Section~\ref{sec:discussion} discusses the results and Section~\ref{sec:conclusion} concludes the paper.

\section{The Current Development Plan}
	\label{sec:nep}

The current NEP in~\cite{Netzentwicklungsplan2017a} targets the year 2030 and comprises several thousand kilometers of new AC and DC transmission lines as well as extensive reinforcement measures, which total to expansion costs of approximately \euro\,$35$\,bn. To study operational properties of the NEP, a detailed system model is required. However, the models employed in~\cite{Netzentwicklungsplan2017a} are published in a limited form only. On that account, this work utilizes the model presented in~\cite{Hewes2016a}. Based on publicly available data, it captures the current German transmission grid as well as the expansion measures in the NEP as accurately as possible and provides a platform for representative studies. The corresponding German transmission grid including the current NEP is shown in Figure~\ref{fig:grid:tso}.

\section{The Hybrid Architecture}
	\label{sec:transition}

\begin{figure}[!t]
\centering
\includegraphics{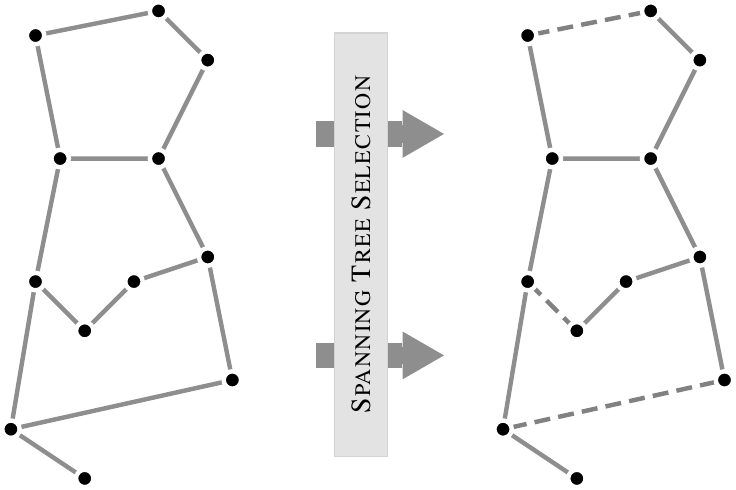}
\caption{The network development strategy suggested by the hybrid architecture comprises the conversion of certain AC branches (solid) to HVDC (dashed). Their selection is based on a spanning tree of the grid's system model.}
\label{fig:transition}
\vskip0.4em
\end{figure}

\begin{figure}[!t]
\centering
\includegraphics{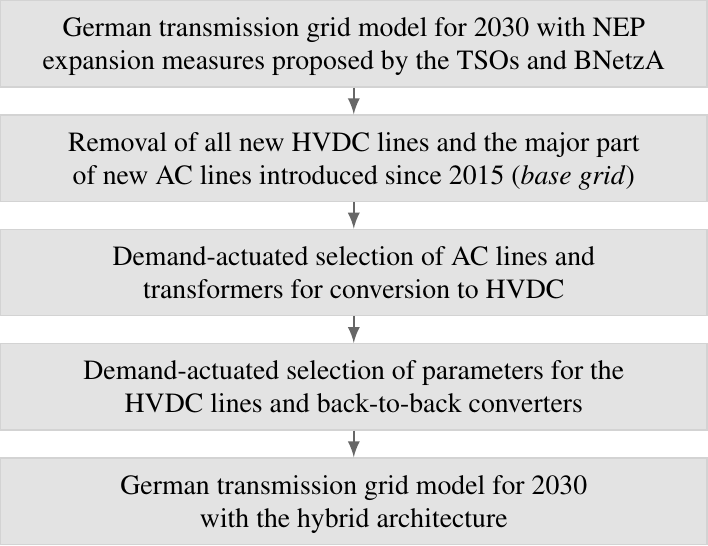}
\caption{Outline of the design procedure of the German transmission grid for 2030 with the hybrid architecture.}
\label{fig:htgdesign}
\end{figure}

In a capacity expansion as proposed by the hybrid architecture, certain AC lines and transformers (branches) are converted to HVDC, where the branches are selected such that loops are resolved~\cite{Hotz2016a}. More formally, the branches that are subject to conversion are selected such that the remaining AC branches constitute a \emph{spanning tree} of the system buses as illustrated in Figure~\ref{fig:transition}. The resulting system structure is referred to as the \emph{hybrid architecture}. Thus, the hybrid architecture proposes a structural transition, which is utilized hereafter to devise an alternative network development strategy. In this process, the NEP is gradually transformed into a transmission grid that features the hybrid architecture, while simultaneously reducing the need for new transmission lines and balancing performance and investment cost, cf. the outline in Figure~\ref{fig:htgdesign}.\vskip-0.4em

\begin{table}[!t]
\centering
\hspace{-1.5mm}\includegraphics{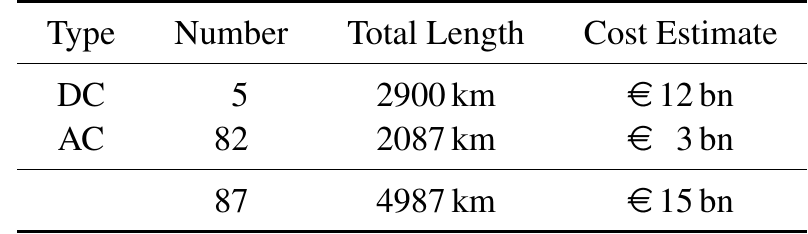}
\caption{Removal of new lines for the base grid.}
\label{tab:base:savings}
\end{table}

\subsection{Base Grid Design}

\begin{figure}[!t]
\centering
\includegraphics{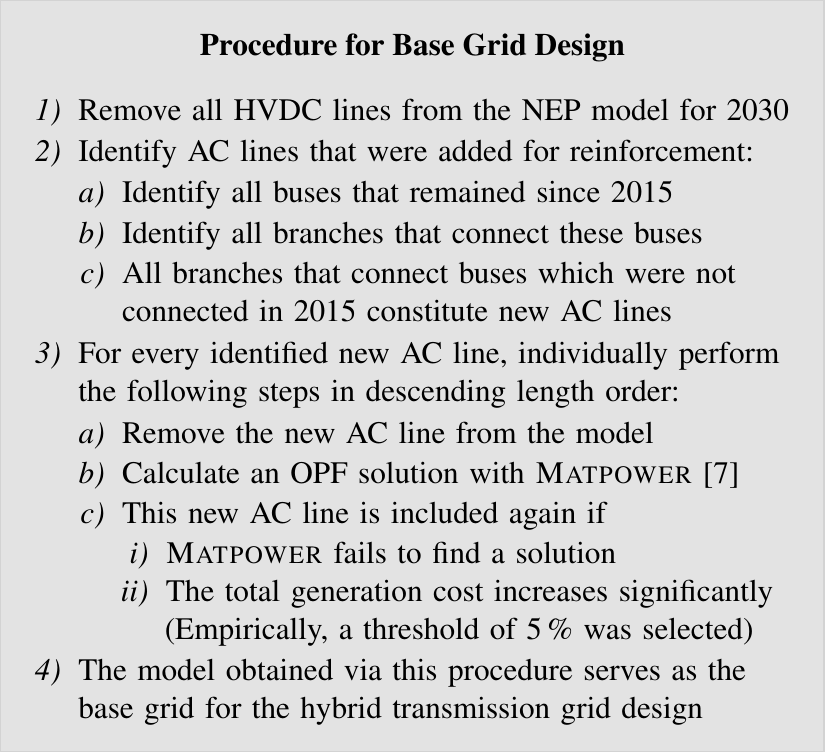}
\vspace{-1.3em}%
\caption{Reduction of the NEP to an appropriate base grid.}
\label{fig:basegridproc}
\end{figure}

To design a transmission grid with the hybrid architecture, an appropriate \emph{base grid} is required that serves as the basis for the structural transition. For the NEP, the transmission grid in 2015 serves as the starting point, where some system buses vanish in 2030 due to decommissioning of certain generation facilities and several new system buses emerge to integrate new generation facilities and reinforce the connection to lower voltage levels. Due to these changes, the grid in 2015 does not constitute an adequate basis for the structural transition, as some expansion measures are essential for a proper integration of 2030's generation and load. For this reason, the NEP is reduced to a base grid as described in Figure~\ref{fig:basegridproc}. This base grid retains only those new lines that are essential to the integration of generation and load. Compared to the NEP, this leads to a substantial reduction of new transmission lines as documented in Table~\ref{tab:base:savings}, where only $6$ out of $93$ new lines are retained.\footnote{\label{foot:corridors}New lines may reside in new or existing corridors. Due to limited information on individual expansion measures in the NEP, the categorization of new AC lines into those in new and existing corridors was not possible. Regarding new DC lines, the NEP considers that $2600$\,km are in new and $300$\,km in existing corridors~\cite[Sec.~4.2.6]{Netzentwicklungsplan2017a}.}

\begin{figure}[!t]
\centering
\includegraphics{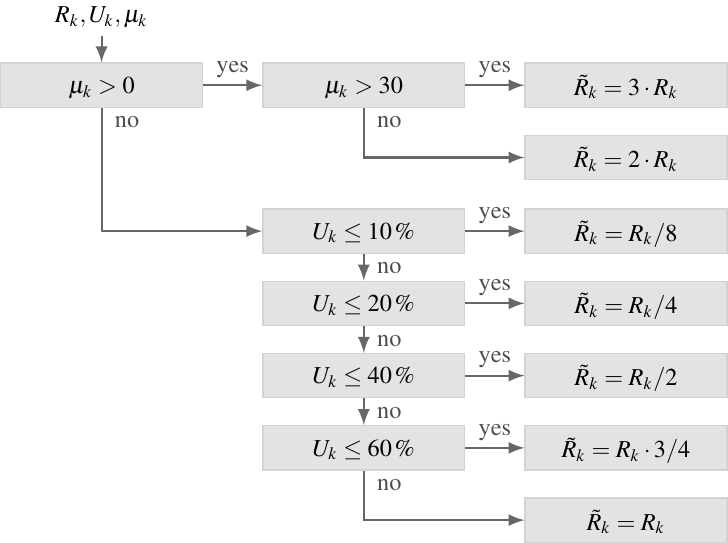}
\caption{The targeted rating $\fTargetRating$ for the conversion of AC branch $k$ is determined by this demand-actuated decision scheme, which is based on the AC branch rating $\fBranchRating$, its utilization $\fBranchUtilization$, and its flow constraint KKT multiplier $\fBranchKktMultiplier$.}
\label{fig:selection:rating}
\end{figure}

\subsubsection{Model Preprocessing}

Prior to the structural transition of the base grid, the model is streamlined via the following preprocessing steps.
\begin{enumerate}[label={\emph{\arabic*)}}, itemsep=-0.15em]
\item Parallel AC branches are combined into an equivalent single AC branch.
\item Multiple generators at a bus are aggregated into an equivalent generator (with a piecewise linear cost).
\item To improve OPF conditioning, AC branches with a series resistance less than $10^{-5}\,$p.u. are imposed with negligible losses by setting their series res.~to~$10^{-5}\,$p.u..
\end{enumerate}

\subsection{Transition to the Hybrid Architecture}

The number of possible transitions to the hybrid architecture equals the number of spanning trees~\cite{Hotz2016a}, which is larger than $10^{308}$ for the base grid. As an exhaustive evaluation of all transitions is intractable, the approach proposed in~\cite[Sec.~VIII-B]{Hotz2017a} is pursued. Therein, every AC branch of the base grid model is associated with an \emph{upgrade suitability measure} to obtain a weighted graph, for which the minimum spanning tree is identified and all AC branches outside this tree are converted to DC operation. To arrive at a demand-actuated upgrade suitability measure $\omega_k$ for AC branch $k$, an OPF is performed for the base grid under peak load using \textsc{Matpower}~\cite{Zimmerman2016a} to obtain the Karush-Kuhn-Tucker (KKT) multiplier $\fBranchKktMultiplier$ of the flow constraint, which is positive if the branch is congested and zero otherwise. This multiplier is augmented by the series resistance $r_k$ of the branch to capture that DC operation becomes more favorable as the line length increases. Finally, it was found experimentally that the cost-performance trade-off is improved if transformers are preferred for the upgrade, which is achieved by augmenting their upgrade suitability measure by the maximum series resistance $r_\mathrm{max}$ of all AC branches. Summarized, the upgrade suitability measure of AC branch $k$ is set to
\begin{equation*}
	\omega_k =
		\begin{cases}
			\fBranchKktMultiplier + r_k
				& \text{if branch $k$ is a line} \\
			\fBranchKktMultiplier + r_k + r_\mathrm{max}
				& \text{if branch $k$ is a transformer\,.}
		\end{cases}
\end{equation*}
\begin{figure}[!t]
\centering
\includegraphics{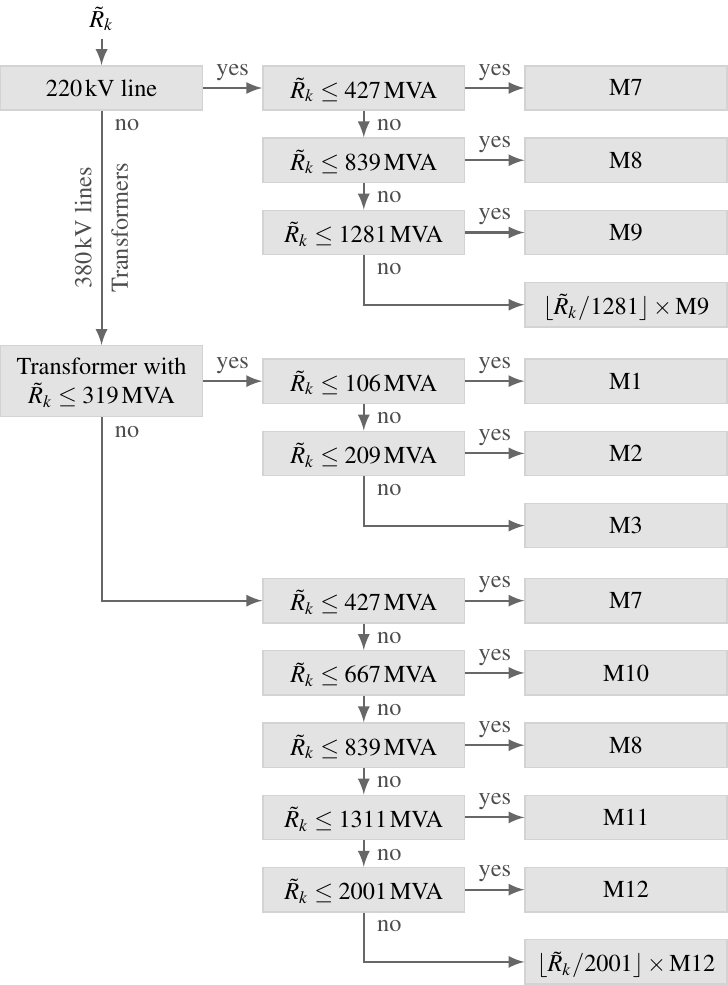}
\caption{For the conversion of AC branch $k$ with target rating $\fTargetRating$, the converter is selected by this decision scheme.}
\label{fig:selection:converter}
\end{figure}
\begin{table}[!t]
\centering
\hspace{-1.5mm}\includegraphics{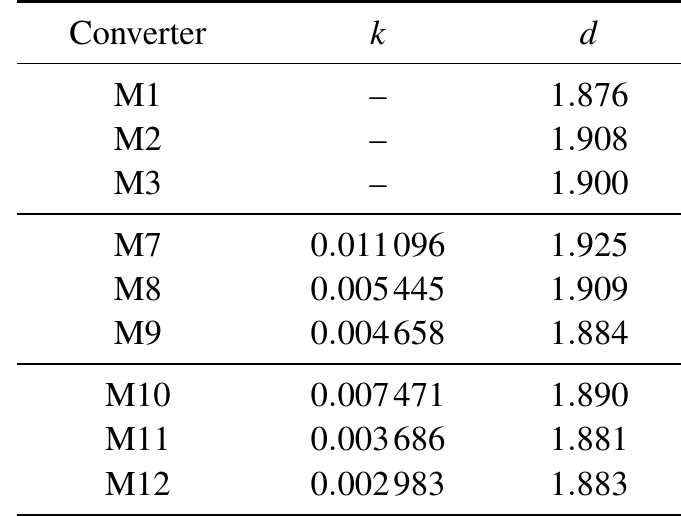}
\caption{Least squares fit of a line with slope $k$ and offset $d$ to the electrical losses with respect to line length to the data in~\cite{ABB-Grid-Systems2012a} for some converters, cf. Fig.~\ref{fig:loss}. For an HVDC line of length $l$, the system model considers losses of $k\cdot l + d$ percent of the active power flow. The converters M1, M2, and M3 are only employed in a back-to-back configuration ($l=0$), where $d$ constitutes the respective loss factor.}
\label{tab:loss}
\end{table}
\begin{figure}[!t]
\centering
\subfloat[\,$\pm$320\,kV symmetric base modules.]{%
\hspace{0.65mm}\includegraphics{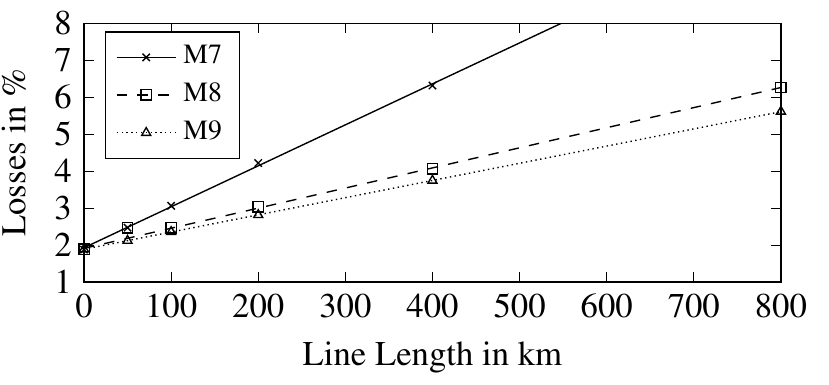}%
\label{fig:loss:320}}\\
\subfloat[\,$\pm$500\,kV symmetric base modules.]{%
\hspace{0.65mm}\includegraphics{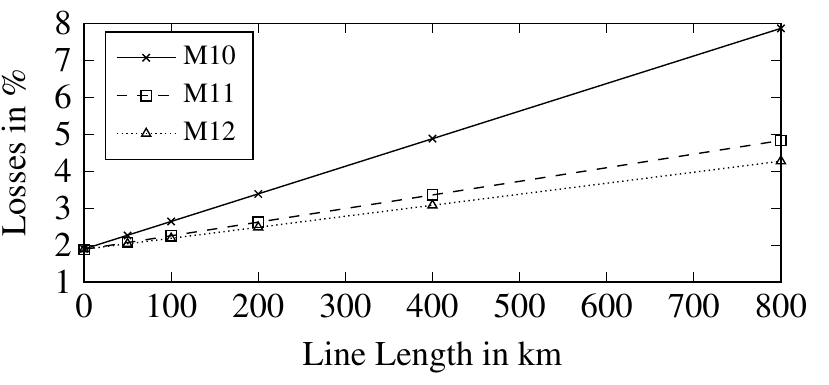}%
\label{fig:loss:500}}%
\caption{Losses of some HVDC Light$^\text{\tiny\textregistered}$ transmission systems by ABB in percent of the active power flow for a typical conductor at base power~\cite{ABB-Grid-Systems2012a}. The markers show the data in~\cite{ABB-Grid-Systems2012a}, while the lines depict the fitted curves in Table~\ref{tab:loss}.}
\label{fig:loss}
\end{figure}
With the minimum spanning tree, $5.98$\,\% of the lines (tot. $4339$\,km) and $55.85$\,\% of the transformers are selected for conversion to an HVDC line and back-to-back (B2B) converter, respectively. To arrive at a demand-actuated target rating $\fTargetRating$ for the conversion of AC branch $k$, its MVA rating $\fBranchRating$ is adapted depending on the severity of congestion, as indicated by the KKT multiplier $\fBranchKktMultiplier$, as well as its utilization $\fBranchUtilization$, which is the percentage of the line flow with respect to the line rating for optimal power flow in the base grid under peak load. The scheme to set the target rating $\fTargetRating$ is shown in Figure~\ref{fig:selection:rating}. Given $\fTargetRating$, an appropriate VSC is selected via the scheme in Figure~\ref{fig:selection:converter}, which exemplarily utilizes converters of ABB~\cite{ABB-Grid-Systems2012a}. In the model, losses of an HVDC line or B2B converter are considered in percent of its active power flow, where the percentage is set according to the converter type and line length based on data by ABB as documented in Figure~\ref{fig:loss} and Table~\ref{tab:loss}.\footnote{The data point at 50\,km for converter M8 in~\cite[p.~29]{ABB-Grid-Systems2012a} appears to be an outlier (typo) and is omitted in the least squares fit in Table~\ref{tab:loss}. The converter M3 exhibits unusually low losses of $0.922\,\%$ in B2B operation~\cite[p.~27]{ABB-Grid-Systems2012a}. This appears to be an irregularity and its losses are assumed~as~$1.9\,\%$.} The model adopts the rating of the converter for the active power limit and considers a reactive power capability of $50\,$\% of this limit, cf.~\cite[Ch.~4]{ABB-Grid-Systems2012a}.

\begin{figure}[!t]
\subfloat[Total load (solid) and RES capacity (dashed).]{%
\hspace{-1mm}\includegraphics{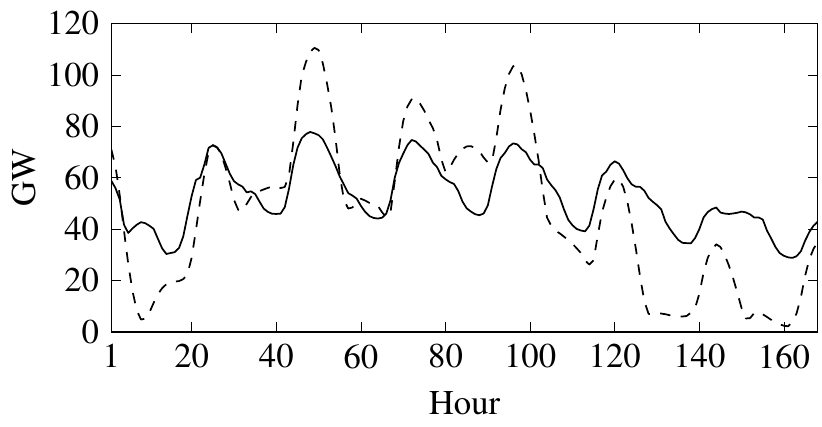}
\label{fig:sim:week2:loadres}}\\[-0.5em]
\subfloat[Total generation cost of the NEP (gray) and HTG (black).]{%
\hspace{-1mm}\includegraphics{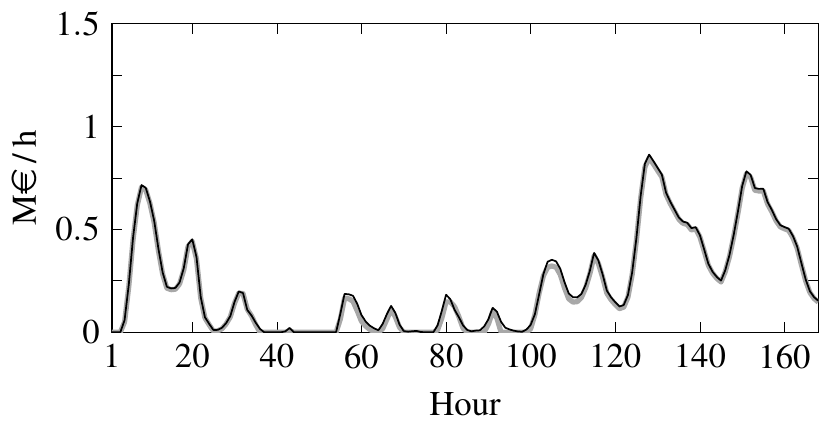}
\label{fig:sim:week2:cost}}%
\caption{Results for the week with year's minimum load.}
\label{fig:sim:week2}
\vskip0.5em
\end{figure}

\begin{figure}[!t]
\subfloat[Voltage profile at peak load of the NEP.]{%
\hspace{-1mm}\includegraphics{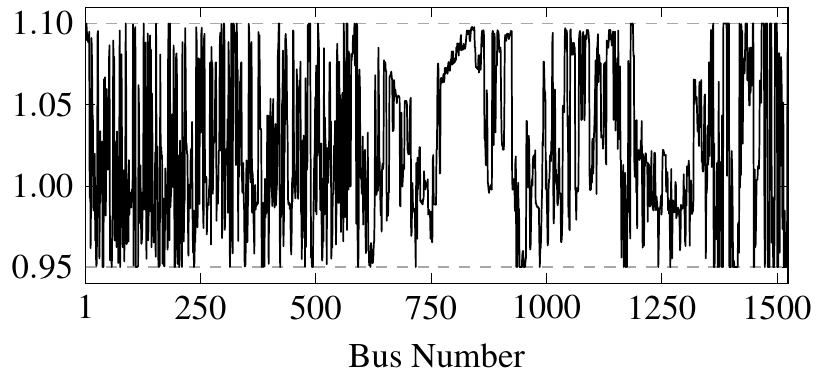}
\label{fig:sim:voltage:max:tso}}\\[-0.25em]
\subfloat[Voltage profile at peak load of the HTG.]{%
\hspace{-1mm}\includegraphics{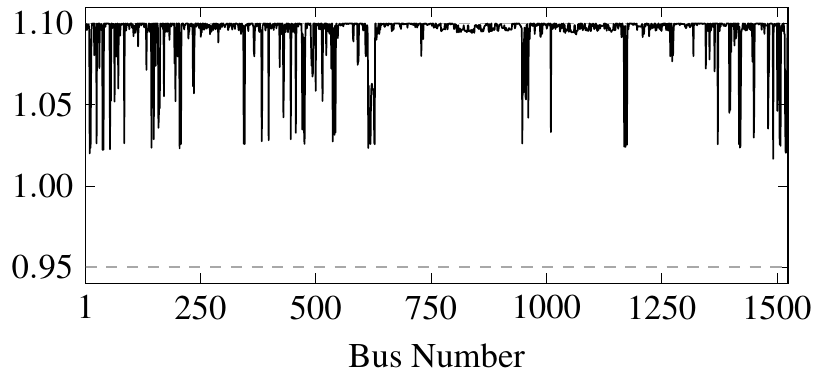}
\label{fig:sim:voltage:max:htg}}%
\caption{Voltage profile at peak load (hour $37$ in Fig.~\ref{fig:sim:week1}). The dashed lines depict the voltage upper and lower bound.}
\label{fig:sim:voltage:max}
\end{figure}

\begin{figure}[!t]
\subfloat[Total load (solid) and RES capacity (dashed).]{%
\hspace{-1mm}\includegraphics{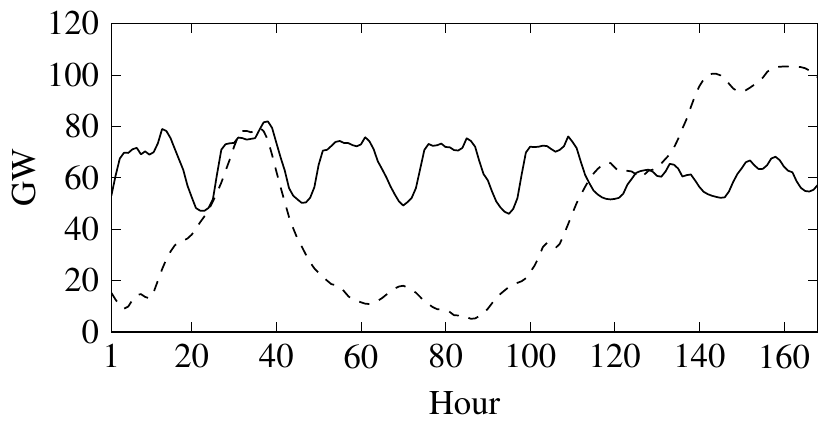}
\label{fig:sim:week1:loadres}}\\[-0.5em]
\subfloat[Total generation cost of the NEP (gray) and HTG (black).]{%
\hspace{-1mm}\includegraphics{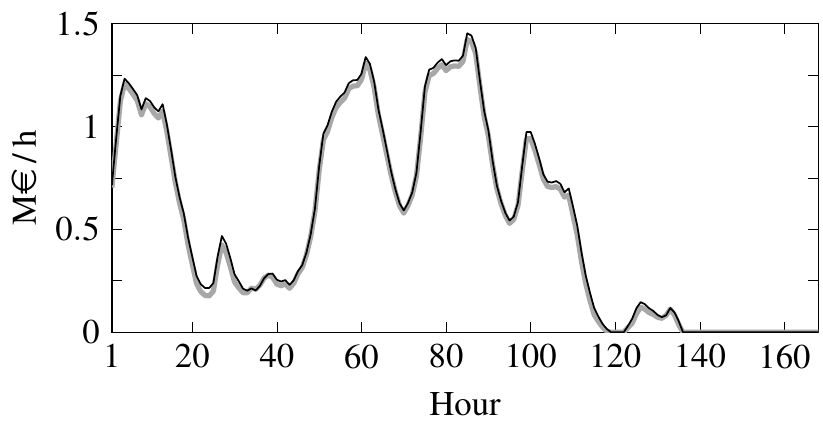}
\label{fig:sim:week1:cost}}%
\caption{Results for the week with year's maximum load.}
\label{fig:sim:week1}
\end{figure}

\begin{table}[!t]
\centering
\hspace{-1.5mm}\includegraphics{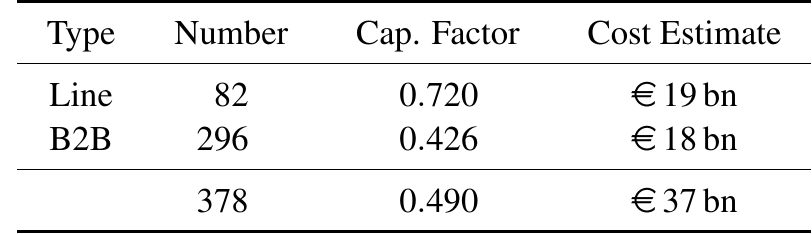}
\caption{Summary of the conversion of AC lines and transformers to HVDC lines and B2B converters. ``Cap. Factor'' denotes the average capacity factor, i.e., the capacity after conversion over the capacity before conversion.}
\label{tab:htg:cost}
\end{table}

\subsection{Performance}

The preceding design process results in an alternative network development strategy for Germany, i.e., a hybrid transmission grid (HTG) that exhibits the \emph{hybrid architecture} as illustrated in Figure~\ref{fig:grid:htg}. To compare its performance against the NEP, an OPF is calculated using \textsc{Matpower}~\cite{Zimmerman2016a} for every hour in the week of minimum and maximum load in the projected year 2030. The results in Figure~\ref{fig:sim:week2} and~\ref{fig:sim:week1} show that both development strategies provide practically the same performance with respect to total generation cost. Due to converter losses, the total generation cost is marginally higher in the HTG, but this is potentially negligible as actual electricity prices consist primarily of charges, taxes, and levies.

\section{Discussion}
	\label{sec:discussion}

The most remarkable feature of the proposed alternative network development strategy is its focus on \emph{existing} corridors, which results in the avoidance of approximately $4900$\,km of new lines compared to the NEP, see also Table~\ref{tab:base:savings}. However, due to the extensive conversion measures and expensive HVDC technology, the transition to the hybrid architecture involves a high cost. For an estimation of the investment volume, the cost assumptions in~\cite{Netzentwicklungsplan2017c} are utilized, i.e., $1.5$\,M\euro/km for new AC lines\footnote{New AC lines in new and existing corridors cannot be distinguished (cf. Footnote~\ref{foot:corridors}) and are considered with the same cost.} as well as new DC overhead lines, $4.0$\,M\euro/km for new DC cables, and $0.2$\,M\euro/km for the conversion of an AC line to DC operation. For HVDC substations, the assumption in~\cite{Netzentwicklungsplan2017c} seems overly conservative, particularly for a large-scale deployment of HVDC. Due to this, the cost estimate of $0.102$\,M\euro/MVA for VSC HVDC substations in~\cite[Sec.~4.2]{Cigre2012a} is adopted. The four new DC lines in the HTG (approx. $440$\,km) are considered as overhead lines, as the NEP also includes overhead lines in the same corridors. On this basis, the savings with respect to the NEP ($100$\,\% cabling of new HVDC lines) and the investment cost for the HTG can be estimated as shown in Table~\ref{tab:base:savings} and Table~\ref{tab:htg:cost}, respectively, which predicts a cost premium of approximately \euro\,$22$\,bn. Concerning this premium, it shall be pointed out that this is a preliminary work and that further considerations may relativize the additional~cost.
\begin{enumerate}[label={\emph{\arabic*)}}, itemsep=-0.15em]
\item For the targeted performance, the branch rating is \emph{reduced} by $50$\,\% on average during conversion to HVDC as shown in Table~\ref{tab:htg:cost}. With a higher target rating during conversion, the hybrid architecture can potentially offer even more transmission capacity on this grid topology. For the year~2035, the TSOs anticipate the necessity of additional HVDC and AC lines~\cite{Netzentwicklungsplan2017a}. To some extent, they may be avoided by the hybrid architecture, which potentially reduces the cost premium.
\item The presented grid design utilizes all expansion measures of the NEP that do not extend the existing grid topology, e.g., reinforcement of existing lines and reactive power compensation. In fact, potentially not all of these measures are necessary, which facilitates a cost reduction. For example, Figure~\ref{fig:sim:voltage:max:htg} suggests that additional compensation is not necessary due to the reactive power provided by the VSCs.
\item The presented grid design is devised with general decision rules that apply equally to all converted branches. With a detailed study and adjustment of individual expansion measures, the investment cost may be reduced.
\item This study employed exclusively point-to-point HVDC systems. If multi-terminal HVDC systems are also considered, the number of VSCs and, thus, the investment cost may be reduced significantly.
\item In the NEP, HVDC lines are implemented as cables (``Erdkabelvorrang''), which results in additional costs of \euro\,$6.5$\,bn compared to overhead lines. As cables entail no technical advantage, this premium is accepted to foster public acceptance. Thus, there is a profound willingness in energy policy to accept extensive investments to support public acceptance, which may even reach out to this alternative expansion strategy.
\end{enumerate}
The hybrid architecture is an ambitious goal, but the preservation of landscape, the flexibility induced by the pervasive incorporation of VSC HVDC systems, and the benefits in operation methods render it highly promising. Due to the extensive conversion measures and their implications, e.g., on grid control, the transition to the hybrid architecture must probably be approached gradually. In this regard, it shall be noted that only $11$ lines and $24$ transformers are uprated during conversion, while all other conversions reduce capacity. Thus, only a minority of the conversions introduces capacity, while the majority only provides controllability and reactive power capacity. Considering that, for this majority, the existing AC lines and transformers offer more capacity, some conversions may not be necessary immediately as the capacity may compensate for reduced controllability. On this basis, it appears promising that appropriate gradual transitions to the hybrid architecture can be devised.

\section{Conclusion}
	\label{sec:conclusion}

In this work, a network development strategy for Germany was presented that focuses on existing corridors and increases capacity by converting systematically selected AC lines and transformers to HVDC. This can avoid up to $4900$\,km of new lines compared to the current NEP of the TSOs and BNetzA and, therewith, reduce the impact on landscape and potentially foster public acceptance. An OPF study of two weeks in the projected year 2030 showed that the performance with respect to total generation cost is on par with the current NEP. The investment cost is considerably higher as HVDC technology is still rather expensive, but its technical properties and alignment with current energy policy may outweigh the cost premium.

The results in this work constitute a preliminary investigation of the potential of the hybrid architecture. Extensive further research is necessary to substantiate its adequacy as a network development strategy. In addition to the points raised in the discussion, a study of its impact on resilience as well as the system's dynamics and stability is of major importance. Regarding the former, the utilization of the controllability of HVDC systems for corrective actions under contingencies is probably essential, necessitating corresponding methods for N-1 secure operational planning. Regarding the latter, the fast and independent control of active and reactive power at this extensive amount of VSC HVDC systems may be coordinated and utilized, e.g., to provide pervasive dynamic voltage support and oscillation~damping.

\setlength{\bibsep}{0pt plus 0.3ex}
\def\bibfont{\small}
\bibliographystyle{unsrt}
\bibliography{htg_germany}

\end{document}